\documentclass{article}%
\usepackage{amssymb}
\usepackage{amsfonts}
\usepackage{graphicx}
\usepackage{amsmath}%
\providecommand{\U}[1]{\protect\rule{.1in}{.1in}}

\newenvironment{proof}[1][Proof]{\textbf{#1.} }{\ \rule{0.5em}{0.5em}}
\begin{document}

\begin{center}
{\LARGE Duality for compact group actions on operator algebras and
applications: irreducible inclusions and Galois correspondence}

\bigskip

\bigskip

Costel Peligrad
\end{center}

\bigskip

\ \ \ \ \ \ \ \ \ \ \ \ \ \ \ \ \ \ \ \ \ \textit{Dedicated to the memory of
Ciprian Foias}

\begin{center}
\bigskip
\end{center}

\bigskip Department of Mathematical Sciences, University of Cincinnati, PO Box
210025, Cincinnati, OH 45221-0025, USA. E-mail address: costel.peligrad@uc.edu

\bigskip


\bigskip


\bigskip

\textbf{ABSTRACT. }We consider compact group actions on C*- and W*- algebras.
We prove results that relate the duality property of the action (as defined in
the Introduction) with other relevant properties of the system such as the
relative commutant of the fixed point algebras being trivial (called the
irreducibility of the inclusion) and also to the Galois correspondence between
invariant C*-subalgebras containing the fixed point algebra and the class of
closed normal subgroups of the compact group.

\begin{center}
\bigskip

\end{center}

\section{Introduction and auxiliary results}

\bigskip

A C*-dynamical system is a triple $(A,G,\alpha)$\ where $A$\ is a C*-algebra,
$G$\ a locally compact group and $\alpha$\ an action of $G$\ on $A,$\ that is,
for each $g\in G,$\ $\alpha_{g}$\ is an automorphism of $A$\ such that the
mapping $g\rightarrow\alpha_{g}(a)$ is continuous for every $a\in A,$ in
particular, $g\rightarrow\varphi(\alpha_{g}(a))$\ is continuous for every
$a\in A$\ and every $\varphi$\ in the dual $A^{\ast}$\ of the C*-algebra
$A.$\ A W*-dynamical system is a triple $(A,G,\alpha)$\ where $A$\ is a
W*-algebra, $G$\ a locally compact group and $\alpha$\ an action of $G$\ on
$A,$ such that the mapping $g\rightarrow\varphi(\alpha_{g}(a))$ is continuous
for every $a\in A$ and every $\varphi\in A_{\ast}\ $where $A_{\ast}\ $is the
predual of $A.$\ In the rest of this paper, we will assume that $G$\ is a
compact group.

\bigskip

In [16], Takesaki proved the following

\bigskip

\textbf{Theorem A }([16], Theorem 1)\textbf{. }Let $(A,G,\alpha)$ be a
W*-dynamical system with $G$\ compact. Suppose that the commutant $S$ of
$\left\{  \alpha_{g}:g\in G\right\}  $ in the group $Aut(A)$\ of all
automorphisms of $A$\ acts ergodically on $A.$ Then, if $\beta\in
Aut(A)$\ leaves $A^{\alpha}$\ pointwise invariant and commutes with $S,$\ it
is of the form $\beta=\alpha_{g}$\ for some $g\in G.$

\bigskip

He also proved a result about a Galois correspondence between $G$\ and
$S$\ invariant subalgebras of $A$ containing $A^{\alpha}$\ and the set of
normal subgroups of $G$ [16, Theorem 3],\ due to Kishimoto [9, Theorem 7]. In
[8] this result is extended to the compact Kac algebras actions on von Neumann
algebras. In Section 3 we will study the Galois correspondence for
C*-dynamical systems.

\bigskip

Notice that if $(A,G,\alpha)$ is a W*-dynamical system with $G$\ compact such
that $(A^{\alpha})^{\prime}\cap A=%
\mathbb{C}
I,$ so if the inclusion $A^{\alpha}\subset A$\ is an irreducible inclusion of
factors, then the automorphism group $S=\left\{  adu:u\in A^{\alpha}\text{
unitary}\right\}  $ acts ergodically on $A$\ and commutes with $\alpha.$
Therefore, a consequence of Theorem A. is the following:

\bigskip

\textbf{Corollary B }[16, Corollary 2]\textbf{. }If the relative commutant
$(A^{\alpha})^{\prime}\cap A$\ of $A^{\alpha}$\ in $A$\ reduces to scalars,
then every automorphism $\beta$ of $A$\ that leaves $A^{\alpha}$\ pointwise
invariant is of the form $\beta=\alpha_{g}$ for some $g\in G.$

\bigskip

The conclusion of Theorem A and Corollary B will be refered in the next
sections as the duality property of the action $\alpha$.

\bigskip

In Section 2. we will prove a converse of Corollary B for W*-dynamical systems
and study the corresponding problem for C*-dynamical systems.

\bigskip

For C*-dynamical systems there are several notions of ergodicity. We mention
some of them that will be used in this paper.

Let $S$\ be a group of automorphisms of a C* algebra $A$\ 

a) $S$\ is called minimal, [10], if the only non zero -invariant hereditary
subalgebra of \ is $A.$

b) $S$\ is said to be topologically transitive, [10], if for every
$B_{1},B_{2}$\ non zero $S$-invariant hereditary subalgebras of $A,$\ we have
$B_{1}B_{2}\neq\left\{  0\right\}  .$ In [2] it is noticed that $S$\ is
topologically transitive if and only if, for every $x,y$ non zero elements of
$A$\ there exists $s\in S$\ such that $xs(y)\neq0.$

Clearly a)$\Rightarrow$b) but not the other way around. The above two notions
generalize the classical concepts for commutative C*-dynamical systems.

In [2] the authors introduced the notion of strong topological transitivity:

c) $S$\ is said to be strongly topologically transitive if for every finite
set $\left\{  x_{i},y_{i}\right\}  _{i=1}^{n}$\ of non zero elements of
$A$\ such that if\ $\sum x_{i}\otimes y_{i}\neq0,$ then there exists $s\in
S$\ such that $\sum x_{i}s(y_{i})\neq0.$ The authors proved a version of
Theorem A for C*-dynamical systems under the strong topological transitivity assumption.

Clearly, c)$\Rightarrow$b) in general, but it is not known whether
b)$\Rightarrow$c) in general. Also we do not know whether a)$\Rightarrow$c) in general.

In order to set the framework for C*-dynamical systems, notice that if
$(A,G,\alpha)$ is a W*-dynamical system with $G$\ compact and $\alpha
$\ faithful such that $(A^{\alpha})^{\prime}\cap A=%
\mathbb{C}
I,$ then every $\alpha_{g},g\neq e$\ is properly outer, i.e. it is not
implemented by a unitary element $u\in A$ (otherwise, $u\in(A^{\alpha
})^{\prime}\cap A=%
\mathbb{C}
I,$ contradiction). For C*-algebras, the proper outerness of automorphisms is
defined [11], [3], using the algebra of local multipliers, $M_{loc}(A),$ (as
denoted in [1]), or equivalently, the algebra of essential multipliers,
$M^{\infty}(A),$ first defined by Pedersen in [11]. This algebra is defined as
the inductive limit of the multiplier algebras $M(I)$\ where $I$ is an
essential ideal of $A$ (an ideal is said to be essential if its annihilator in
$A$\ equals $\left\{  0\right\}  $)$.$\ In this definition, it is assumed, as
proven in [12, Proposition 3. 12. 8.], that if $I_{1}\subset I_{2}$ are
essential ideals then there is a unit preserving inclusion $M(I_{2})\subset
M(I_{1}).$

\bigskip

In the next lemma, if $B$\ is a C*-algebra, $B^{\ast\ast}$\ stands for the
second dual of $B$. For each automorphism $\alpha$\ of $B,$\ we denote by
$\alpha^{\ast\ast}$\ the second dual of $\alpha.$

\bigskip

\textbf{1.1. Lemma }\textit{Let }$(B,G,\alpha)$\textit{\ be a C*-dynamical
system and }$M(B)$\textit{\ the multiplier algebra of }$B.$\textit{\ Then,}

\textit{i) For every }$x\in M(B)$\textit{\ and }$\varphi\in B^{\ast}%
$\textit{\ the mapping }$g\rightarrow\varphi(\alpha_{g}^{\ast\ast}%
(x))$\textit{\ is continuous.}

\textit{ii) For every }$x\in M(B)$ \textit{there exists a unique element
}$P^{\alpha}(x)\in M(B)^{\alpha}$\textit{ such that}%
\[
\varphi(P^{\alpha}(x))=\int\varphi(\alpha_{g}^{\ast\ast}(x))dg
\]
\textit{for\ every }$\varphi\in B^{\ast}.$

\textit{iii) The mapping }$P^{\alpha}:M(B)\rightarrow M(B)^{\alpha}$
\textit{is a (norm one) conditional expectation continuous in the }$B^{\ast}$
\textit{topology of }$M(B)$\textit{.}

\textit{iv) If }$\rho$ \textit{is a faithful representation of }%
$B$\textit{\ then the normal extension }$\rho"$\textit{ of }$\rho$\textit{\ to
}$B^{\ast\ast}$\textit{\ is faithful on }$M(B),$\textit{\ the function
}$g\rightarrow\varphi(\rho"(\alpha_{g}^{\ast\ast}(x)))$\textit{\ is continuous
for every }$x\in M(B)$\textit{\ and }$\varphi\in(\rho"(B^{\ast\ast}))_{\ast}%
$\textit{\ where }$(\rho"(B^{\ast\ast}))_{\ast}$\textit{\ denotes the predual
of }$\rho"(B^{\ast\ast})$\textit{ and}%
\[
\varphi(\rho^{"}(P^{\alpha}(x)))=\int\varphi(\rho"(\alpha_{g}^{\ast\ast
}(x)))dg
\]
\textit{for every }$x\in M(B)$\textit{ and }$\varphi\in(\rho"(B^{\ast\ast
}))_{\ast}.$

\bigskip

\begin{proof}
i) Let $\widetilde{B}$ be C*-algebra generated by $B$\ and the unit of
$B^{\ast\ast}$ and $\widetilde{B}_{sa}$\ the self adjoint part of
$\widetilde{B}.$ Then, by [12, 3.12.9.] we have
\[
M(B)_{sa}=(\widetilde{B}_{sa})^{m}\cap(\widetilde{B}_{sa})_{m}%
\]
where $(\widetilde{B}_{sa})^{m}$ (respectively $(\widetilde{B}_{sa})_{m})$
denotes the strong limits in $B^{\ast\ast}$\ of all increasing nets in
$\widetilde{B}_{sa}$ (respectively the limits of all decreasing nets in
$\widetilde{B}_{sa}).$ Let $x\in M(B)_{sa}.$ Then there exists an increasing
net $\left\{  \widetilde{a}_{\nu}=a_{\nu}+c_{\nu}I\right\}  _{\nu}$ and a
decreasing net $\left\{  \widetilde{b}_{\mu}=b_{\mu}+d_{\mu}I\right\}  ,$ both
in $\widetilde{B}_{sa}$ that converge strongly to $x$\ in $B^{\ast\ast}.$
Clearly, for every $\nu$ and $\mu$ and for every $\varphi\in B^{\ast},$ the
mappings $g\rightarrow\varphi(\alpha_{g}^{\ast\ast}(\widetilde{a}_{\nu
}))=\varphi(\alpha_{g}(a_{\nu}))+c_{\nu}\varphi(I)$\ and $g\rightarrow
\varphi(\alpha_{g}^{\ast\ast}(\widetilde{b}_{\mu}))=\varphi(\alpha_{g}(b_{\mu
}))+d_{\mu}\varphi(I)$\ are continuous and $\lim_{\nu}\varphi(\alpha_{g}%
^{\ast\ast}(\widetilde{a}_{\nu}))=\lim_{\mu}\varphi(\alpha_{g}^{\ast\ast
}(\widetilde{b}_{\mu}))=\varphi(\alpha_{g}^{\ast\ast}(x)).$ Therefore, if
$\varphi\in B^{\ast}$\ is a positive functional, the mapping $g\rightarrow
\varphi(\alpha_{g}^{\ast\ast}(x))$\ is both lower and upper semi continuous
and therefore continuous. Since this is true for positive functionals, it
follows that it holds for every $\varphi\in B^{\ast}.$

ii) In [5, Theorem 1.1.] it is proven that both dual pairs of Banach spaces
($M(B),B^{\ast})$ and ($B^{\ast},M(B))$ posess the Krein property and
therefore the corresponding (Pettis) integral can be performed. They also
proved that the $B^{\ast}$\ topology and the strict topology of $M(B)$\ are consistent.

iii) This follows from i) and ii).

iv) The fact that $\rho"$ is faithful on $M(B)$\ follows from [12, Corollary
3.12.5.]. To prove the continuity of the mapping \ and the subsequent
equality, notice that if $\varphi\in(\rho"(B^{\ast\ast}))_{\ast}$, then the
composition $\varphi\circ\rho"\in B^{\ast}$ and we can apply i).
\end{proof}

\bigskip

Now, let $(B,G,\alpha)$\ be a C*-dynamical system with $G$\ compact and $B$\ a
prime C*-algebra (i.e. a C*-algebra, $B$ such that every ideal of $B$\ is
essential). In [3, the proof of implication 6$\Rightarrow$7 and\ Proposition
3.3.] it is proven that every ideal of $B$\ contains an $\alpha$-invariant
ideal and therefore, $M_{loc}(B)$\ is the inductive limit of the multiplier
algebras $M(I)$\ where $I$ is a non zero $\alpha$-invariant ideal of $B.$ In
what follows we will assume that $B$\ has a faithful factorial
representation.It is known that this is the case if, in particular, $B$\ is a
separable prime C*-algebra. We can assume that $B\subset B(H)$\ for some
Hilbert space, $H,$\ and $\mathcal{M}=B"$ is a factor, where $B"$\ is the
bicommutant of $B$. If this is the case, then for every ideal $I\subset B$\ we
have $I"=\mathcal{M}.$ Therefore the predual $\mathcal{M}_{\ast}$ of
$\mathcal{M}$\ is canonically embedded in the dual $I^{\ast}$\ of $I.$ Also,
according to [12, Corollary 3.12.5], the multiplier algebra $M(I)$\ of
$I$\ can be canonically embedded in $\mathcal{M}$ for every ideal $I\subset
B$\ and, therefore, $M_{loc}(B)$\ can be identified with the norm closure
\[
\overline{\cup\left\{  M(I):I\subset B\text{ \textit{is an }}\alpha
\text{\textit{-invariant ideal}}\right\}  }\subset\mathcal{M}.
\]
Since $B"$\ is a factor, and therefore for every non zero ideal $I\subset
B,$\ $I"=B",$ if $(B")_{\ast}$ denotes the predual of $B",$\ it follows that
$(B")_{\ast}\subset I^{\ast}$\ for every non zero ideal $I$\ of $B.$ Since by
Lemma 1.1. iv) the (faithful) representation\ restricted to $I$\ can be
extended to a normal representation of $I^{\ast\ast}$\ which is faithful on
$M(I),$\ for every non zero ideal $I\subset B$\ and $x\in M(I)\subset B",$\ we
can identify $\alpha_{g}^{\ast\ast}(x)\in I^{\ast\ast}$ with$\ $an element
$\alpha_{g}^{"}(x)\in B".$\ Using these facts and Lemma 1.1. for $B=I,$ where
$I\subset B$ is an $\alpha$-invariant ideal, we have:

\bigskip

\textbf{1.2. Lemma }\textit{For every }$x\in M_{loc}(B)$\textit{\ there exists
a unique }$P^{\alpha}(x)\in M_{loc}(B)$\textit{\ such that the mapping
}$g\rightarrow\varphi(\alpha_{g}^{"}(x))$\textit{\ is continuous for every
}$\varphi\in(B")_{\ast}$\textit{ and}%
\[
\varphi(P^{\alpha}(x))=\int\varphi(\alpha_{g}^{"}(x))dg
\]
\textit{for all }$\varphi\in(B")_{\ast}$\textit{ and }$P^{\alpha}$\textit{\ is
a conditional expectation of }$M_{loc}(B)$\textit{\ onto the fixed point
algebra }$(M_{loc}(B))^{\alpha}.$

\bigskip

\begin{proof}
This follows from Lemma 1.1. iv) and the fact that if $I\subset B$ is an
$\alpha$-invariant ideal, we have $M(I)^{\alpha}=M(I^{\alpha}).$ This latter
equality follows from the fact that every approximate identity of $I^{\alpha}%
$\ is also an approximate identity of $I.$
\end{proof}

\bigskip

\textbf{1.3.} \textbf{Lemma }\textit{Let }$B$\textit{\ be a C*-algebra.}
\textit{If }$x\in M_{loc}(B)$\textit{\ is such that }$BxB=\left\{  0\right\}
,$ \textit{then }$\mathit{x=0.}$

\bigskip

\begin{proof}
By definition, $M_{loc}(B)=\overline{\cup_{I}M_{I}}$ where, for each essential
ideal $I\subset B,$ $M_{I}$\ is a unital C*-algebra isomorphic with the
multiplier algebra $M(I)$\ of $I.$\ Then, if $x\in M_{loc}(B)$, there is a
sequence $\left\{  x_{n}\right\}  \subset\cup_{I}M_{I},x_{n}\in M_{I_{n}}$
such that $lim_{n}x_{n}=x$ in norm. If $x\neq0$\ we can assume that
$\left\Vert x\right\Vert =1$\ and $\left\Vert x_{n}\right\Vert =1$\ for all
$n\in%
\mathbb{N}
.$ Let $N\in%
\mathbb{N}
$\ be such that $\left\Vert x_{n}-x\right\Vert <\frac{1}{2}$ for all
$n\geqslant N.$\ Then, since by hypothesis, $BxB=\left\{  0\right\}  ,$\ we
have $ax_{n}b=a(x_{n}-x)b$\ for every $a,b\in B$ and every $n,$\ in particular
for $n=N$\ and every $a,b\in I_{N}$ such that $\left\Vert a\right\Vert
\leq1,\left\Vert b\right\Vert \leq1.$ Therefore, $\left\Vert ax_{N}%
b\right\Vert <\frac{1}{2}$ for all $a,b\in I_{N}$ with $\left\Vert
a\right\Vert \leq1,\left\Vert b\right\Vert \leq1.$ However, since $x_{N}\in
M_{I_{N}}$\ and $M_{I_{N}}$ is isomorphic with $M(I_{N}),$ we have
$\sup\left\{  \left\Vert ax_{N}b\right\Vert :a,b\in I_{N},\left\Vert
a\right\Vert \leq1,\left\Vert b\right\Vert \leq1\right\}  =\left\Vert
x_{N}\right\Vert =1$ (this can be seen by taking for $a,b$ two distinct
elements of an approximate identity of $I_{N}$ and then taking the
supremum$).$ This contradiction shows that $x=0.$
\end{proof}

\bigskip

For simplicity, in the next section we will write $\alpha_{g}(x)$\ instead of
$\alpha_{g}^{"}(x)$\ for $x\in M_{loc}(B).$\ Also, the element $P^{\alpha
}(x),x\in M_{loc}(B),$\ from Lemma 1.2. will be denoted, simply%
\[
P^{\alpha}(x)=\int\alpha_{g}(x)dg.
\]

If $G$\ is a compact group, the set of equivalence classes of unitary
irreducible representations of $G$ is denoted by $\widehat{G}$ and called dual
object (or dual space) of the group $G.$ For each $\pi\in\widehat{G}$\ we
denote also by $\pi$\ a representative of that class. If $\pi$\ is an
irreducible unitary representation of $G,$\ and $d_{\pi}$ the dimension of the
corresponding Hilbert space, $H_{\pi},$\ $\chi_{\pi}(g)=d_{\pi}^{-1}\sum
_{i=1}^{d_{\pi}}\overline{\pi_{ii}(g)}$ denotes its character.

Now let $(A,G,\alpha)$ be a W*- or a C*-dynamical system.

\bigskip

\textbf{1.4}. \textbf{Notations }\textit{i) For each }$\pi$\textit{\ denote}%
\[
A_{1}^{\alpha}(\pi)=\left\{  \int\chi_{\pi}(g)\alpha_{g}(a)dg:a\in A\right\}
\subset A
\]
\textit{the corresponding spectral subspace of }$A.$\textit{\ }

\textit{ii) }$A_{2}^{\alpha}(\pi)=\left\{  \left[  a_{ij}\right]  \in A\otimes
B(H_{\pi}):(\alpha_{g}\otimes\iota)(\left[  a_{ij}\right]  )=\left[
a_{ij}\right]  (1\otimes\pi_{g}),g\in G\right\}  .$

\textit{In [13, JFA] it is noticed that each entry, }$a_{ij},$\textit{ is the
element of }$A_{1}^{\alpha}(\pi)$\textit{ such that }$\alpha_{g}(a_{ij}%
)=\sum_{l}\pi_{lj}(g)a_{il}.$

\textit{iii) }$\overline{\pi}$\textit{ denotes the conjugate representation of
}$\pi.$\textit{\ Clearly, }$A_{1}^{\alpha}(\overline{\pi})=A_{1}^{\alpha}%
(\pi)^{\ast}.$

\textit{iv) The Arveson spectrum spectrum of the action is defined as
}$sp(\alpha)=\left\{  \pi\in\widehat{G}:A_{1}^{\alpha}(\pi)\neq\left\{
0\right\}  \right\}  .$

\textit{v) The linear span of }$\left\{  m^{\ast}n:m,n\in A_{2}^{\alpha}%
(\pi)\right\}  $\textit{ is denoted by }$A_{2}^{\alpha}(\pi)^{\ast}%
A_{2}^{\alpha}(\pi)$\textit{ and the linear span of }$\left\{  mn^{\ast
}:m,n\in A_{2}^{\alpha}(\pi)\right\}  $\textit{ is denoted }$A_{2}^{\alpha
}(\pi)A_{2}^{\alpha}(\pi)^{\ast}.$\textit{ It is straightforward to check that
}$A_{2}^{\alpha}(\pi)A_{2}^{\alpha}(\pi)^{\ast}$\textit{ is a two sided ideal
of the fixed point algebra }$(A\otimes B(H_{\pi}))^{\alpha\otimes\iota}%
$\textit{ and }$A_{2}^{\alpha}(\pi)^{\ast}A_{2}^{\alpha}(\pi)$\textit{ an
ideal of }$(A\otimes B(H_{\pi}))^{\alpha\otimes ad\pi}.$ \textit{If }%
$A_{2}^{\alpha}(\pi)^{\ast}A_{2}^{\alpha}(\pi)$\textit{\ is dense in
}$(A\otimes B(H_{\pi}))^{\alpha\otimes ad\pi}$\textit{ for all }$\pi
\in\widehat{G},$\textit{\ the action }$\alpha$\textit{\ is called saturated.}

\bigskip

\section{Duality for compact group actions versus irreducible inclusions.}

\bigskip

In this section we will state and prove a converse of Corollary B for W*-- and
C*-dynamical systems and also prove the Corollary B for some C*-dynamical
systems. The following example shows that a verbatim converse of Corollary B
is false.

\bigskip

\textbf{2.1. Example }Let $K$ be a finite dimensional Hilbert space of
dimension larger than $1$ and $G$ the compact group of all unitary operators
in $A=B(K).$ If $\alpha_{g}(x)=gxg^{-1},g\in G,x\in A,$\ then $A^{\alpha}=%
\mathbb{C}
I,$ so $\left(  A^{\alpha}\right)  ^{\prime}\cap A=A\neq%
\mathbb{C}
I,$\ but the action has the duality property.

\bigskip

In the above example, clearly, both $A^{\alpha}$ and $A$\ are factors so even
if we add the condition that $A^{\alpha}$ and $A$\ are factors to the duality
property the conclusion that $\left(  A^{\alpha}\right)  ^{\prime}\cap A=%
\mathbb{C}
I$ is false for non abelian compact groups. We mention that if $G$\ is a
compact abelian group then the duality property plus the assumption that
$A^{\alpha}$\ and $A$\ are factors imply that the inclusion $A^{\alpha}\subset
A$\ is irreducible, i.e. $\left(  A^{\alpha}\right)  ^{\prime}\cap A=%
\mathbb{C}
I.$ This fact was proven in [3,Theorem]. Thus we will be concerned with
compact non abelian groups $G.$ We start with the following

\bigskip

\textbf{2.2. Lemma }\textit{Let }$(A,G,\alpha)$\ \textit{be a C*-dynamical
system (respectively a W*-dynamical system) with }$G$\textit{\ compact. If
(}$A\otimes B(H_{\pi}))^{\alpha\otimes ad\pi}$ \textit{are prime C*-algebras
(respectively W* factors) for all }$\pi\in sp(\alpha),$\textit{ then
}$sp(\alpha)$\textit{ is closed under tensor products. If in addition,
}$\alpha$\textit{\ is faithful (that is }$\alpha_{g}\neq\iota$ if $g\neq e),$
\textit{then} $sp(\alpha)=\widehat{G}.$

\bigskip

\begin{proof}
Recall that, by definition, $sp(\alpha)$ is said to be closed under tensor
products if every irreducible component of the tensor product $\pi_{1}%
\otimes\pi_{2},\pi_{1},\pi_{2}\in sp(\alpha)$\ is in $sp(\alpha).$ Since
\textit{(}$A\otimes B(H_{\pi}))^{\alpha\otimes ad\pi}$\ are prime
(respectively factors), it follows that $A_{2}^{\alpha}(\pi)^{\ast}%
A_{2}^{\alpha}(\pi),\pi\in sp(\alpha),$\ are essential ideals of
\textit{(}$A\otimes B(H_{\pi}))^{\alpha\otimes ad\pi}$ (respectively weakly
dense ideals in the case of W*-dynamical systems).\ Thus, $sp(\alpha
)$\ coincides with the spectrum defined in [6, Definition 3.1. (1)] for
compact quantum groups. Applying [6, Lemma 5.5] to the particular case of
groups, it follows that $sp(\alpha)$ is closed under tensor products. Since
$sp(\alpha)$\ is also closed under conjugation by the definition of the
spectrum (since $A_{1}(\overline{\pi})=A_{1}(\pi)^{\ast}),$ we can apply [7,
Theorem 28.9.] and derive that the annihilator $sp(\alpha)_{\bot}$\ of
$sp(\alpha)$\ is a closed normal subgroup $G_{0}$\ of $G$\ and the annihilator
$G_{0}^{\perp}$\ of $G_{0}$\ in $\widehat{G}$\ equals $sp(\alpha).$ It follows
that, if $\alpha$\ is faithful, $G_{0}=\left\{  e\right\}  ,$\ so
$sp(\alpha)=\widehat{G}.$
\end{proof}

\bigskip

As we mentioned in Section 1, in this Section 2, if $(A,G,\alpha)$\ is a
C*-dynamical system with $G$\ compact, we will assume that $A$\ has a faithful
factorial representation, or, equivalently, as we stated before, that
$A\subset B(H)$\ for some Hilbert space $H$\ and $A"$\ is a factor, so, in
particular, $A$\ is a prime C*-algebra. This assumption will allow us to use
the identification of $M_{loc}(A)$\ with a norm closed subalgebra of $A^{"}%
$\ and use Lemma 1.2. for $B=A$. Notice that if $(A,G,\alpha)$\ is a
W*-dynamical system, we do not need the algebra $M_{loc}(A)$\ and thus, in
Corollary 2.4. below we will assume that $A$\ is a factor. We will also assume
that the action $\alpha$\ is faithful.

The next result is our converse of Corollary B for C*-dynamical systems. Our
version of Corollary B itself will be given in Theorem 2.5.

\bigskip

\textbf{2.3. Theorem} \textit{Let }$(A,G,\alpha)$\textit{ be a C*-dynamical
system with }$G$\textit{\ compact, }$\alpha$ \textit{faithful and }%
$A$\textit{\ as above. If}

\textit{i) }$(A\otimes B(H_{\pi}))^{\alpha\otimes ad\pi},\pi\in sp(\alpha
)$\textit{\ are prime C*-algebras and}

\textit{ii) If }$\beta$\textit{\ is an automorphism of }$M_{loc}%
(A)$\textit{\ which leaves }$A^{\alpha}$\textit{\ pointwise invariant\ then
}$\beta=\alpha_{g}$\textit{\ for some }$g\in G,$

\textit{then }$sp(\alpha)=\widehat{G}$\textit{ and }$(A^{\alpha})^{\prime}\cap
M_{loc}(A)=%
\mathbb{C}
I.$

\bigskip

\begin{proof}
The equality $sp(\alpha)=\widehat{G}$\ follows from the above Lemma 2.2. We
will prove next the second part of the conclusion. Let $u\in(A^{\alpha
})^{\prime}\cap M_{loc}(A)$\ be a unitary element. Then $\beta
(a)=adu(a)=uau^{\ast}$ is an automorphism of $M_{loc}(A)$ that leaves
$A^{\alpha}$ pointwise invariant. By the hypothesis ii) there exists $g_{0}\in
G$\ such that $\beta=\alpha_{g_{0}}.$ Let $\pi\in\widehat{G}$\ and $m\in
A_{2}^{\alpha}(\pi),m\neq0.$ Then, by the definition of $A_{2}^{\alpha}(\pi
)$\ we have ($\alpha_{g}\otimes\iota)(m)=m(I\otimes\pi_{g})$\ for every $g\in
G,$ in particular, for $g=g_{0}.$\ Since $\alpha_{g_{0}}=\beta=adu$\ it
follows that
\[
(u\otimes1)m=m(u\otimes\pi_{g_{0}}).\text{
\ \ \ \ \ \ \ \ \ \ \ \ \ \ \ \ \ \ (1)}%
\]
Clearly, $u\otimes\pi_{g_{0}}\in M_{loc}(A)\otimes B(H_{\pi})=M_{loc}(A\otimes
B(H_{\pi})).$\ This latter equality follows from the fact that every ideal of
$A\otimes B(H_{\pi})$\ is of the form $I\otimes B(H_{\pi})$\ where $I$\ is an
ideal of $A.$ Let $n\in A_{2}^{\alpha}(\pi).$\ By multiplying \ (1) on the
right by $n^{\ast}$\ we get%
\[
(u\otimes1)mn^{\ast}=m(u\otimes\pi_{g_{0}})n^{\ast}.\text{
\ \ \ \ \ \ \ \ \ \ \ \ \ \ \ (2)}%
\]
By applying $\alpha_{g}\otimes\iota,g\in G$\ to equality (2), since $mn^{\ast
}\in A^{\alpha}\otimes B(H)$\ and then using the definition of $A_{2}^{\alpha
}(\pi)$\ in the right hand side of (2), we get%
\[
(\alpha_{g}(u)\otimes1)mn^{\ast}=m((\alpha_{g}\otimes ad\pi_{g})(u\otimes
\pi_{g_{0}})n^{\ast}.\text{ \ (3)}%
\]
By integrating both sides of (3), we get%
\[
(P^{\alpha}(u)\otimes I)mn^{\ast}=mP^{\alpha\otimes ad\pi}(u\otimes\pi_{g_{0}%
})n^{\ast}.\text{ \ \ \ \ \ \ (4)}%
\]
Since $u\in(A^{\alpha})^{\prime}\cap M_{loc}(A)$\ and this latter algebra is
$\alpha$-invariant\ we have that $P^{\alpha}(u)\in(A^{\alpha})^{\prime}\cap
M_{loc}(A)^{\alpha}.$\ On the other hand, according to [3, Proposition 3.3.],
$M_{loc}(A)^{\alpha}\subset M_{loc}(A^{\alpha}),$ so $P^{\alpha}%
(u)\in(A^{\alpha})^{\prime}\cap M_{loc}(A^{\alpha}).$ Since by hypothesis
$A^{\alpha}$\ is a prime C*-algebra, from [3, Proposition 3.1.] it follows
that $P^{\alpha}(u)=\lambda I$ for some scalar $\lambda.$ Therefore%
\[
\lambda(I\otimes I)mn^{\ast}=mP^{\alpha\otimes ad\pi}(u\otimes\pi_{g_{0}%
})n^{\ast}.\text{ \ \ \ \ \ \ \ \ (4}^{,}\text{)}%
\]
Let $c,d\in A_{2}^{\alpha}(\pi).$\ By multiplying (4') on the right by
$c$\ and on the left by $d^{\ast}$\ we have%
\[
\lambda(I\otimes I)d^{\ast}mn^{\ast}c=d^{\ast}mP^{\alpha\otimes ad\pi
}(u\otimes\pi_{g_{0}})n^{\ast}c.\text{ \ \ \ \ \ \ \ \ (5)}%
\]
We noticed in Section 1 that $J_{\pi}=\overline{A_{2}^{\alpha}(\pi)^{\ast
}A_{2}^{\alpha}(\pi)}$ is an ideal of $(A\otimes B(H_{\pi}))^{\alpha\otimes
ad\pi}.$ Therefore, from (5) it follows that%
\[
J_{\pi}(\lambda(I\otimes I)-P^{\alpha\otimes ad\pi}(u\otimes\pi_{g_{0}%
}))J_{\pi}=\left\{  0\right\}  .\text{ \ \ \ \ \ \ \ \ \ \ \ \ (6)}%
\]
So, if we denote $x=\lambda(I\otimes I)-P^{\alpha\otimes ad\pi}(u\otimes
\pi_{g_{0}})\in(M_{loc}(A\otimes B(H_{\pi}))^{\alpha\otimes ad\pi}$ we have%
\[
J_{\pi}xJ_{\pi}=\left\{  0\right\}  .\text{
\ \ \ \ \ \ \ \ \ \ \ \ \ \ \ \ \ \ \ \ \ \ \ \ \ \ \ \ \ \ \ \ \ \ \ \ \ \ \ \ \ \ \ \ \ \ \ \ \ \ \ (6')}%
\]
and since by [3, Proposition 3.3.], $M_{loc}(A\otimes B(H_{\pi})^{\alpha
\otimes ad\pi}\subset M_{loc}((A\otimes B(H_{\pi}))^{\alpha\otimes ad\pi}),$
and $(A\otimes B(H_{\pi}))^{\alpha\otimes ad\pi}$\ is a prime C*-algebra by
hypothesis i), it follows that $J_{\pi}\ $is an essential ideal of $(A\otimes
B(H_{\pi}))^{\alpha\otimes ad\pi},$ so, by [11], $M_{loc}((A\otimes B(H_{\pi
}))^{\alpha\otimes ad\pi}=M_{loc}(J_{\pi}).$\ By taking $B=J_{\pi}$ in Lemma
1.3., it follows that $x=0.$ Therefore%
\[
\lambda(I\otimes I)=\int(\alpha_{g}(u)\otimes\pi_{gg_{0}g^{-1}})dg.\text{
\ \ \ \ \ \ \ \ \ \ \ \ \ \ \ \ \ \ \ \ \ \ (7)}%
\]
Let us calculate the sum of diagonal elements in (7):%
\[
\lambda d_{\pi}=\sum_{l,k,j}\int\pi_{lk}(g)\pi_{kj}(g_{0})\pi_{jl}%
(g^{-1})\alpha_{g}(u)dg=
\]%
\[
=\sum_{k,j}\pi_{kj}(g_{0})\int\sum_{l}\pi_{jl}(g^{-1})\pi_{lk}(g)\alpha
_{g}(u)dg=
\]%
\[
\sum_{k,j}\pi_{kj}(g_{0})\int\delta_{kj}\alpha_{g}(u)dg
\]
where $\delta_{kj}$\ is the Kronecker symbol. Therefore%
\[
\lambda d_{\pi}=tr(\pi_{g_{0}})P^{\alpha}(u)=\lambda tr(\pi_{g_{0}}).\text{
\ \ \ \ \ \ \ \ \ \ \ \ \ \ \ \ \ \ \ \ \ \ (8)}%
\]
Hence, if $\lambda\neq0$\ then $tr(\pi_{g_{0}})=d_{\pi}$ for all $\pi
\in\widehat{G}.$ It follows that $\pi_{g_{0}}=I$\ for every $\pi\in\widehat
{G}.$ Since $\widehat{G}$ separates the points of $G$\ we get $g_{0}=e,$ the
neutral element of $G.$ Therefore, $\beta=adu=\iota$, so $u\in A^{\prime}\cap
M_{loc}(A).$\ Since $A$\ is prime by hypothesis, from [3, Proposition 3.1.] it
follows that $u$\ is a scalar multiple of the identity. Now, let
$a\in(A^{\alpha})^{\prime}\cap M_{loc}(A)$\ be a positive element of norm $1.$
Then, there exists a unitary element, $u\in(A^{\alpha})^{\prime}\cap
M_{loc}(A)$ such that $a=\frac{1}{2}(u+u^{\ast}).$\ Since $P^{\alpha}$\ is
faithful, it follows that $P^{\alpha}(a)\neq0$\ and therefore, $P^{\alpha
}(u)\neq0.$ Acording to the previous arguments, $u$\ and therefore, $a$\ is a
scalar multiple of the identity and the proof is completed.
\end{proof}

\bigskip

The following is a consequence of the proof of the above Proposition.

\bigskip

\textbf{2.4. Corollary }\textit{Let }$(A,G,\alpha)$\textit{ be a W*-dynamical
system with }$G$\textit{\ compact. If}

\textit{i) }$(A\otimes B(H_{\pi}))^{\alpha\otimes ad\pi},\pi\in sp(\alpha
)$\textit{\ and }$A$\textit{\ are factors and}

\textit{ii) If }$\beta$\textit{\ is an automorphism of }$A$\textit{\ which
leaves }$A^{\alpha}$\textit{\ pointwise invariant\ then }$\beta=\alpha_{g}%
$\textit{\ for some }$g\in G,$

\textit{then }$sp(\alpha)=\widehat{G}$\textit{ and }$(A^{\alpha})^{\prime}\cap
A=%
\mathbb{C}
I.$

\bigskip

\begin{proof}
If $(A\otimes B(H_{\pi}))^{\alpha\otimes ad\pi},\pi\in sp(\alpha)$ are
factors, then $A_{2}^{\alpha}(\pi)^{\ast}A_{2}^{\alpha}(\pi),\pi\in
sp(\alpha),$\ are weakly dense ideals of \textit{(}$A\otimes B(H_{\pi
}))^{\alpha\otimes ad\pi}$ and the rest of the proof of Theorem 2.3. carries
over to the case of W*-dynamical systems.
\end{proof}

\bigskip

Next we will prove that hypothesis of the above corollary is equivalent with
the irreducibility of the inclusion $A^{\alpha}\subset A$\ for W*-dynamical systems.

\bigskip

\textbf{2.5. Theorem }\textit{Let }$(A,G,\alpha)$\textit{ be a W*-dynamical
system with }$G$\textit{\ compact and }$\alpha$\textit{\ faithful. The
following statements, i) and ii), are equivalent}

\textit{i) (i1) }$(A\otimes B(H_{\pi}))^{\alpha\otimes ad\pi},\pi\in
sp(\alpha)$\textit{\ and }$A$\textit{\ are factors.}

\textit{ \ \ (i2) If }$\beta$\textit{\ is an automorphism of }$A$%
\textit{\ which leaves }$A^{\alpha}$\textit{\ pointwise invariant\ then
}$\beta=\alpha_{g}$\textit{\ for some }$g\in G.$

\textit{ii) }$(A^{\alpha})^{\prime}\cap A=%
\mathbb{C}
I.$

\textit{Under each of the conditions i) and ii) we have} $sp(\alpha
)=\widehat{G}.$

\bigskip

\begin{proof}
The implication i)$\Rightarrow$ii) and the fact that i) implies the equality
$sp(\alpha)=\widehat{G}$\ follows from the above Corollary. Suppose ii) holds.
Then, according to [17, Proposition 6.2.], applied to the particular case of
compact groups, the relative commutant of $A$\ in the crossed product
$A\times_{\alpha}G$\ equals $%
\mathbb{C}
I.$ In particular, this means that $A$ and $A\times_{\alpha}G$\ are factors.
Applying [13 Cor. 3.12.] to the case of W*-dynamical systems (as noticed at
the end of the paper), it follows that $sp(\alpha)=\widehat{G}\ $and
$(A\otimes B(H_{\pi}))^{\alpha\otimes ad\pi},\pi\in\widehat{G}$ are factors,
so i1) holds.\ The condition i2) follows from Corollary B of Takesaki, stated
in Section 1.
\end{proof}

\bigskip

The next result is concerned with C*-dynamical systems. To prove the converse
of Theorem 2.3., i.e. a version of Corollary B for C*-dynamical systems, we
need a condition that is stronger than the existence of a faithful factorial
representation of $A.$ Namely, we will assume that there exists a faithful
representation $\rho$\ of $A$\ such that $\rho(A^{\alpha})^{\prime}\cap
\rho(A)"=%
\mathbb{C}
I.$\ Clearly such a representation is factorial and, since $M_{loc}%
(A)\subset\rho(A)",$ it implies that $(A^{\alpha})^{\prime}\cap M_{loc}(A)=%
\mathbb{C}
I$. Examples of such systems are given in [3] for compact abelian groups
$G$\ and a stronger condition than this, namely the existence of a faithful
irreducible representation $\rho$\ of $A$\ such that the restriction of $\rho
$\ to $A^{\alpha}$\ is also irreducible, for non abelian compact groups in
[4]. So, we will assume that $A\subset B(H)$\ for some Hilbert space $H$\ and
$(A^{\alpha})^{\prime}\cap A^{"}=%
\mathbb{C}
I.$

\bigskip

\textbf{2.6. Proposition} \textit{Let }$(A,G,\alpha)$\textit{ be a
C*-dynamical system with }$G$\textit{\ compact and }$\alpha$%
\textit{\ faithful. Suppose that }$A\subset B(H)$\ for some Hilbert space
$H$\ and $(A^{\alpha})^{\prime}\cap A^{"}=%
\mathbb{C}
I..$\textit{Then,}

\textit{i) (i1)\ }$sp(\alpha)=\widehat{G}$\textit{, }$(A\otimes B(H_{\pi
}))^{\alpha\otimes ad\pi},\pi\in\widehat{G}$\textit{\ are prime and}

\textit{ \ (i2) If }$\beta$\textit{\ is an automorphism of }$M_{loc}%
(A)$\textit{\ which leaves }$A^{\alpha}$\textit{\ pointwise invariant\ then
}$\beta=\alpha_{g}$\textit{\ for some }$g\in G.$

\bigskip

\begin{proof}
Applying [3, Corollary 5.4.], it follows that the hypothesis implies that the
unitary group of $A^{\alpha}$ (with unit adjoined if necessary)\ acts strongly
topologically transitively on $A$ (strong topological transitivity as defined
in [2, Introduction]). Since obviously the action of the unitary group of
$A^{\alpha}$\ commutes with $\alpha$ we can apply [2, Corollary 2.2.] to
obtain the equality $sp(\alpha)=\widehat{G}.$ Since $(A^{\alpha})^{\prime}\cap
A^{"}=%
\mathbb{C}
I$ and obviously $(A^{\alpha})"\subset A"$ it follows that $(A^{\alpha})"$\ is
a factor.\ From the discussion in Section 1 it follows that $M_{loc}%
(A^{\alpha})\subset(A^{\alpha})",$\ so, in particular, $(A^{\alpha})^{\prime
}\cap M_{loc}(A^{\alpha})=%
\mathbb{C}
I.$ From [3, Proposition 3.1.] it follows that $A^{\alpha}$\ is a prime
C*-algebra. Since $\alpha$\ commutes with a strongly topologically transitive
action and $A^{\alpha}$\ is prime, we can apply [14, Proposition 5.4.] and we
get that the crossed product $A\times_{\alpha}G$\ is also a prime C*-algebra.
From [13, Corollary 3.12.] it follows that $(A\otimes B(H_{\pi}))^{\alpha
\otimes ad\pi},\pi\in\widehat{G}$\textit{\ }are prime so i1) is proven. i2)
follows from the updated version of [2, Theorem 2.1.] given in [3, proof of
the implication 12$\Rightarrow$13].
\end{proof}

\bigskip

\section{Duality for compact actions and Galois correspondence}

\bigskip

If $(A,G,\alpha)$ is a C*-dynamical system with $G$\ compact, and $B$\ is an
$\alpha$-invariant C* subalgebra of $A$\ such that $A^{\alpha}\subset B\subset
A,$\ in [15] it is defined a closed normal subgroup of $G$%
\[
G^{B}=\left\{  g\in G:\alpha_{g}(b)=b,b\in B\right\}
\]
and a subalgebra of $A$%
\[
A^{G^{B}}=\left\{  a\in A:\alpha_{g}(a)=a,g\in G^{B}\right\}  .
\]

In [15, Lemma 13 and Corollary 9] we proved that if there exists a subgroup
$S\subset Aut(A)$\ which acts minimally on $A,$ that commutes with $\alpha
,$\ then $A^{G^{B}}=B$. We have also shown that this is false if $S$ acts
strongly topologically transitively on $A.$ In what follows we will complete
that result as follows. Suppose as above that there exists a subgroup
$S\subset Aut(A)$\ which acts\ minimally on $A,$ that commutes with $\alpha.$
If $G_{0}\subset G$\ is a closed normal subgroup of $G,$\ let%
\[
A^{G_{0}}=\left\{  a\in A:\alpha_{g}(a)=a,g\in G_{0}\right\}
\]
and%
\[
G^{A^{G_{0}}}=\left\{  g\in G:\alpha_{g}(a)=a,a\in A^{G_{0}}\right\}  .
\]

In Proposition 3.2. below, we will prove that $G^{B}=G_{0}$\ where
$B=A^{G_{0}}.$ This will follow from the following

\bigskip

\textbf{3.1. Theorem} \textit{Let }$(A,G,\alpha)$\textit{ be a C*-dynamical
system with }$G$\textit{\ compact and }$\alpha$\textit{\ faithful.\ Suppose
that there exists a subgroup }$S\subset Aut(A)$\textit{\ which acts\ minimally
on }$A,$\textit{ that commutes with }$\alpha.$\textit{ Suppose as above that
there exists a subgroup }$S\subset Aut(A)$\textit{\ which acts\ minimally on
}$A,$\textit{ and commutes with }$\alpha.$ \textit{Then, if }$\beta
$\textit{\ is an automorphism of }$A$\textit{\ that commutes with }%
$S$\textit{\ and leaves }$A^{\alpha}$\textit{\ pointwise invariant, there
exists }$g\in G$\textit{\ such that }$\beta=\alpha_{g}.$

\bigskip

\begin{proof}
Applying [15, Lemma 13], it follows that, under the hypothesis of the theorem
the action $\alpha$\ is saturated, i.e. $A_{2}^{\alpha}(\pi)^{\ast}%
A_{2}^{\alpha}(\pi)$\ is dense in $(A\otimes B(H_{\pi}))^{\alpha\otimes ad\pi
}$ for all $\pi\in\widehat{G}.$ Since $S\subset Aut(A)$\textit{ }%
acts\ minimally on $A,$ it is in particular topologically transitive and
therefore, by [10, Proposition 2.14.]\ the only fixed points of the extension
of $S$\ to the multiplier algebra $M(A)$\ are scalars multiples of the
identity. In [10], [14] we have refered to this property as weak ergodicity of
$S.$\ The conclusion of the theorem follows from [14, Theorem 3.3.].
\end{proof}

\bigskip

\textbf{3.2. Proposition }\textit{Let }$(A,G,\alpha)$\textit{ be a
C*-dynamical system with }$G$\textit{\ compact and }$\alpha$%
\textit{\ faithful.\ Suppose that there exists a subgroup }$S\subset
Aut(A)$\textit{\ which acts\ minimally on }$A,$\textit{ that commutes with
}$\alpha.$ \textit{Then, if }$G_{0}$\textit{\ is a closed normal subgroup of
}$G$\textit{,} $B=A^{G_{0}}.$\textit{is }$S$\textit{- and }$\alpha
$\textit{-invariant, }$A^{\alpha}\subset B$\textit{ and }$G^{B}=G_{0}$.

\bigskip

\begin{proof}
Clearly $A^{\alpha}\subset B=A^{G_{0}}\subset A\ $and$,$ since $S$\ commutes
with $\alpha,$\ $B$\ is $S$-invariant. Since $G_{0}$\ is a normal subgroup,
$B$\ is also $\alpha$-invariant. From the definitions of $B$\ and $G^{B}$\ we
get that $G_{0}\subset G^{B}.$\ Now consider the C*-dynamical system
$(A,G_{0},\alpha|_{G_{0}})$ and let $g\in G^{B}.$\ Then $\alpha_{g}$ leaves
$A^{G_{0}}=B$\ pointwise invariant and commutes with the minimal action
$S.$\ Applying Theorem 3.1., it follows that there exists $g_{0}\in G_{0}%
$\ such that $\alpha_{g}=\alpha_{g_{0}}.$\ Since $\alpha$\ is faithful,
$g=g_{0}\in G_{0},$ so $G^{B}=G_{0}.$
\end{proof}

\bigskip

The next result completes and improves the results in [15].

\bigskip

\textbf{3.3. Theorem} \textit{Let }$(A,G,\alpha)$\textit{ be a C*-dynamical
system with }$G$\textit{\ compact and }$\alpha$\textit{\ faithful.\ Suppose
that there exists a subgroup }$S\subset Aut(A)$\textit{\ which acts\ minimally
on }$A,$\textit{ that commutes with }$\alpha.$ \textit{Then the assignments
}$B\rightarrow G^{B}$\textit{ and }$G_{0}\rightarrow A^{G_{0}}$\textit{ are
inverse to each other maps between the set of all }$S$\textit{- and }$\alpha
$\textit{-invariant C* subalgebras of }$A$ \textit{that contain }$A^{\alpha}%
$\textit{\ and the set of all closed normal subgroups of }$G.$

\bigskip

\begin{proof}
From [15, Lemma 13 and Corollary 9] it follows that if $B$\ is $S$- and
$\alpha$-invariant and $A^{\alpha}\subset B,$\ then $A^{G^{B}}=B.$\ From
Proposition 3.2. above it follows that if $G_{0}$\ is a closed normal subgroup
of $G,$\ then $G^{A^{G_{0}}}=G_{0}$\ and the proof is completed.
\end{proof}

\bigskip

\textbf{References}

\bigskip

[1] P. Ara and M. Mathieu, \textit{Local multipliers of C*-algebras,
}Springer, 2003\textit{.}

[2] O. Bratteli, G. A. Elliott and D. W. Robinson, \textit{Strong topological
transitivity and C*-dynamical systems, }J. Math. Soc. Japan, 37 (1985), 115-133.

[3] O.Bratteli, G. A. Elliott, D. E. Evans and A. Kishimoto,
\ \textit{Quasi-product actions of a compact abelian group on a C*algebra,}
Tohoku J. of Math., 41 (1989), 133-161.

[4] O.Bratteli, G. A. Elliott and A. Kishimoto, \textit{Quasi-product actions
of a compact groups on a C*algebra,} J. Funct. Anal. 115 (1993), 313-343.

[5] C. D'Antoni and L. Zsido, \textit{Groups of linear isometries on
multiplier C*-algebras,} Pacific J. Math. 193 (2000), 279--306.

[6] R. Dumitru and C. Peligrad, \textit{Spectra for compact quantum group
actions and crossed products, }Transactions of the Amer. Math. Soc., 364
(2012), 3699-3713.

[7] E. Hewitt and K. A. Ross, \textit{Abstract Harmonic Analysis, vol. II,
}Springer, 1970.

[8] M. Izumi, R. Longo and, S. Popa, \textit{A Galois correspondence for
compact group of automorphisms of von Neumann algebras with a generalization
to Kac algebras, }J. Funct. Anal. 155 (1998), 25-63.

[9] A. Kishimoto, \textit{Remarks on compact automorphism groups of a certain
von Neumann algebra, }Publ. RIMS, Kyoto Univ., 13 (1977), 573-581.

[10] R. Longo and C. Peligrad, \textit{Noncommutative topological dynamics and
compact actions on C*-algebras, }J. Funct. Anal.,58 (1984), 157-175.

[11] G. K. Pedersen, \textit{Approximating derivations on ideals of a
C*-algebra, }Inventiones Math. 45 (1978), 299-305.

[12] G. K. Pedersen, \textit{C*-algebras and their automorphism groups,
}Academic Press, 1979.

[13] C. Peligrad, \textit{Locally compact group actions on C*-algebras and
compact subgroups, }J. Funct. Anal. 76 (1988), 126-139.

[14] C. Peligrad, \textit{Compact actions commuting with ergodic actions and
applications to crossed products,} Transactions of the Amer. Math. Soc., 331
(1992), 825-836.

[15] C. Peligrad, \textit{A Galois correspondence for compact group actions on
C*-algebras, }J. Funct. Anal., 261 (2011), 1227-1235.

[16] M. Takesaki, \textit{Fourier analysis of compact automorphism groups (an
application of Tannaka duality theorem), }Coll Internat. du CNRS, 274 (1979).

[17] S. Vaes, \textit{The unitary implementation of a locally compact group
action, }J. Funct. Anal. 180 (2001), 426-480.

\end{document}